\documentclass[12pt]{amsart}
\usepackage{amsmath,amscd,amsthm,amsfonts, amssymb,amsxtra, mathrsfs,enumitem,relsize, stackengine,soul,cancel}
\usepackage{calrsfs}
\usepackage{caption}
\usepackage[us,12hr]{datetime}
\usepackage[all]{xy} \SelectTips{cm}{}
\usepackage{hyperref}
 \hypersetup{colorlinks=true,citecolor=blue}
   \usepackage{graphicx}
     \usepackage[margin=1.5in]{geometry}
   \usepackage[font=scriptsize]{caption}
\usepackage{todonotes}
\usepackage{tikz}
\usetikzlibrary{shapes.geometric,decorations.pathreplacing}
\usetikzlibrary{positioning}
\CDat
\usepackage{scalerel}[2016/12/29]

 \subjclass{55Q25}
\newtheorem{thm}{Theorem}[section]  
\newtheorem*{un-no-thm}{Theorem}

\newtheorem{bigthm}{Theorem}

\theoremstyle{definition}

\theoremstyle{definition}

\theoremstyle{definition}

\theoremstyle{definition}

\theoremstyle{remark}
\newtheorem{rem}[thm]{Remark}

\newtheorem*{acks}{Acknowledgements}

\DeclareMathOperator{\Sp}{Spectra}

\begin{document}
\title{A note on the second James-Hopf invariant}
\date{\today}
\author{John R.\ Klein}
\address{Wayne State University, Detroit, MI 48202}
\email{klein@math.wayne.edu}
\begin{abstract}  This paper characterizes the stabilized second James-Hopf invariant by means of three axioms. Specifically, we show that it is the unique 
natural transformation satisfying the Cartan formula, {\color{blue} and} vanishing on suspensions\st{, and a metastable EHP property}. 
The proof combines the natural stable splitting of the James construction with Goodwillie calculus.
\end{abstract}
\maketitle
\setlength{\parindent}{15pt}
\setlength{\parskip}{1pt plus 0pt minus 1pt}

\def\Sp{\text{\bf Sp}}
\def\vo{\varOmega}
\def\vs{\varSigma}
\def\smsh{\wedge}
\def\flush{\flushpar}
\def\id{\text{id}}
\def\dbslash{/\!\! /}
\def\codim{\text{\rm codim\,}}
\def\:{\colon}
\def\holim{\text{holim\,}}
\def\hocolim{\text{hocolim\,}}
\def\Bbb{\mathbb}
\def\bold{\mathbf}
\def\Aut{\text{\rm Aut}}
\def\cal{\mathcal}
\def\sec{\text{\rm sec}}
\def\gda{G\text{\rm -}\delta\text{\rm -}\alpha}
\def\PDD{\text{\rm pd\,}}
\def\PD{\text{\rm P}}
\def\stableto {\,\, \mapstochar \!\!\to}

\setcounter{tocdepth}{1}


\section{Introduction}

{\color{blue} Nick Kuhn as pointed out to me that the EHP property (iii) from the previous arXiv submission is unnecessary. I am indicating the modifications in blue and striking out text where needed.}
\medskip

The James-Hopf invariants \cite{James-suspension-triad} are natural transformations
\[
\cal J_n \: [\Sigma A,\Sigma B] \to [\Sigma A,\Sigma B^{[n]}]\, ,
\]
in which $A$ and $B$ are based spaces, $B^{[n]}$ denotes the $n$-fold smash product of $B$,
and $[X,Y]$ is the set of homotopy classes of based maps  $X \to Y$. 

Let $E\: [X,Y] \to [\Sigma X,\Sigma Y]$ be the one-fold suspension map.
Boardman and Steer \cite{Boardman-Steer}
consider the one-fold suspended James invariants 
\[
h_n := E\cal J_n\: [\Sigma A,\Sigma B] \to [\Sigma^2 A,\Sigma^2 B^{[n]}]
\]
and show they form a  {\it Hopf ladder.} The latter consists of three axioms:
\begin{enumerate}
\item $h_1$ is the identity;
\item $h_n(Ef) = 0$ for $f\in [A,B]$;
\item The collection $\{h_n\}$ satisfies the {\it Cartan formula}:
\[
h_n (\alpha + \beta) = \sum_{i+j=n} h_i(\alpha)\cdot h_j(\beta)\, ,\qquad \alpha,\beta \in [\Sigma A,\Sigma B]\, .
\]
\end{enumerate}
We note that (ii) is really an infinite collection of axioms
rolled into one. By \cite[thm~2.2]{Boardman-Steer}, there is precisely one Hopf ladder. 

In this paper we will consider the fully stabilized second James-Hopf invariant
\[
\gamma \: [\Sigma A,\Sigma B] \to \{A,B\smsh B\}
 \]
in which $\gamma(f)$ coincides with the stable homotopy class of $\cal J_2(f)$. 
 Then $\gamma (Ef) = 0$ for $f \in [A,B]$ and  the {\it Cartan Formula} in this case is the single statement
 \[
 \gamma(f+g) = \gamma(f) + f\cup g + \gamma(g)\, , \qquad f,g\in [\Sigma A,\Sigma B]\, ,
 \]
 where the {\it cup product} $f\cup g$ is the stable homotopy class of the composition
 \[
 \Sigma^2 A @> \Sigma \Delta_A >> \Sigma^2 A\smsh A  = (\Sigma A) \smsh (\Sigma A) @> f \smsh f >> \Sigma B\smsh \Sigma B = \Sigma^2 B\smsh B\, .
 \]
\st{The invariant $\gamma$ also possesses the \mbox{\it EHP property:} There is a global constant $c\in \Bbb Z$ such that
if}
\begin{enumerate}
\item \st{$A$ is a CW complex of dimension $\le k$,}
\item \st{ $B$ is $r$-connected, }
\item \st{ $k \le 3r+c$, and}
\item \st{$\gamma(f)=0$,}
\end{enumerate}
\st{Then $f = Eg$ for some $g\in [A,B]$. In the case of $\gamma$, we may take $c=1$.}

\begin{bigthm} \label{bigthm:axioms} Suppose that $\lambda\: [\Sigma A,\Sigma B] \to \{A,B\smsh B\}$ is a natural transformation
satisfying
\begin{enumerate}[label=(\roman*).]
\item $\lambda\circ E = 0$, {\color{blue} and}
\item $\lambda$ satisfies the Cartan formula, \st{and}
\item \st{$\lambda$ possesses the EHP property.}
\end{enumerate}
then $\lambda = \gamma$.
\end{bigthm}

\begin{rem}  Even when restricted to the case of $h_2$, Boardman and Steer's proof \cite[thm.~3.15]{Boardman-Steer}
requires  the Cartan formula in the case of $h_n$ for every $n$.
 \st{In Theorem}  \cancel{\ref{bigthm:axioms}}\st{, we are able to avoid the higher Hopf invariants at the cost of requiring the EHP property.}
\end{rem}

\begin{acks} This paper was motivated by insightful comments by a referee  of the paper \cite{klein2026stablehopf}, to whom I am grateful.
{\color{blue} I am  indebted to Nick Kuhn for informing me that the EHP property is unnecessary, and providing the current proof which is very short and clean. See the remark after the proof
for the related approach of \cite[App.~B]{Kuhn_diag}.}
\end{acks}

\section{Proof of Theorem \ref{bigthm:axioms}}
Let $B$ be a based space. Recall that $J(B)$ is a reduced free monoid on the points of $B$.
In \cite{James_reduced_product}, it was proved that the natural map
\[
J(B) \to \Omega \Sigma B
\]
is a homotopy equivalence when $B$ is a connected CW complex. (If $B$ is not connected, then the map
is known to be group completion.)

Moreover, $J(B)$ comes equipped with a filtration 
\[
J_1(B) \subset J_2(B) \subset \cdots
\]
in which $J_k(B)\subset J(B)$ is the subspace defined by 
the words of length $\le k$.
Then
\begin{enumerate}[label=(\roman*).]
\item $J_1(B) = B$, and
\item $J_k(B)/J_{k-1}(B) = B^{[k]}$.
\end{enumerate}
Furthermore, the cofibration sequence of spectra
\[
\Sigma^{\infty} J_{k-1}(B) \to \Sigma^{\infty} J_k(B) \to \Sigma^\infty B^{[k]}
\]
naturally splits (cf.~\cite[ex.~1.20]{Goodwillie_calc3}), so we obtain natural  weak equivalence of spectra
\[
\textstyle\bigvee\limits_{n\ge 1} \Sigma^\infty B^{[n]} \,\, \simeq \,\, \Sigma^\infty \Omega \Sigma B\, .
\]
It follows that the natural transformation $\lambda$ is determined by a natural transformation
\[
\textstyle\bigvee\limits_{n\ge 1} \lambda_{2,n}\: \textstyle\bigvee\limits_{n\ge 1} \Sigma^\infty B^{[n]}\to \Sigma^\infty B^{[2]} \, ,
\]
where $\lambda_{2,n} \: \Sigma^\infty B^{[n]}\to \Sigma^\infty B^{[2]}$. 

The homotopy functor $B\mapsto B^{[n]}$ is homogeneous of degree $n$, and 
by standard Goodwillie calculus arguments \cite{Goodwillie_calc3}, $\lambda_{2,n}$ is homotopically trivial for $n > 2$.
Furthermore,  $\lambda_{2,1}\: \Sigma^\infty B \to \Sigma^\infty B^{[2]} $ is trivial since $\lambda\circ E$ is trivial.

Note that in the case of the stabilized James-Hopf invariants, $\gamma_{2,1}$ is the identity by definition.
Hence, we only need to establish that 
\[
\lambda_{2,2}\:  \Sigma^\infty B^{[2]}\to \Sigma^\infty B^{[2]} 
\]
is homotopically the identity.

{\color{blue} Let $p_1, p_2 \colon B \times B \to B$ be the two
projections. Then the Cartan formula (ii)  implies that $\lambda(p_1 + p_2) = p_1 \cup p_2 \in \{(B\times B)_+,B\smsh B\}$.
That is, it is the homotopy class of the composition
\[
\Sigma^\infty (B\times B)_+ @> q >>  \Sigma^\infty B\smsh B @>\lambda_{2,2} >>  \Sigma^\infty B\smsh B 
\]
where $q\: (B\times B)_+ = B_+ \smsh B_+ \to B\smsh B$ is the projection obtained by mapping $+$ to the basepoint of $B$.
As the induced homomorphism
\[
q^\ast\: \{B\smsh B, B\smsh B\} \to \{(B\times B)_+,B\smsh B\}
\]
is split injective and $q^*(1_{B \smsh B}) = q^\ast(\lambda_{2,2})$, the result is immediate. \qed}

\st{The spectrum $\Sigma^\infty B^{[2]}$ is homogenous of degree two. Its derivative is the 
naive $\Bbb Z_2$-spectrum $S[\Bbb Z_2] := \Sigma^{\infty}{\Bbb Z_2}_+$.  Since $\lambda_{2,2}$ is a natural transformation of homotopy functors,
it is classified up to homotopy by the homotopy class of a $\Bbb Z_2$-equivariant self-map $S[\Bbb Z_2]\to S[\Bbb Z_2]$.
The ring of homotopy classes of such maps is to the group ring $\Bbb Z[\Bbb Z_2]$.
For $\theta\in \Bbb Z[\Bbb Z_2]$, we will write}
\[
\cancel{\hat \theta\: S[\Bbb Z_2]\to S[\Bbb Z_2]} 
\] 
\st{for the corresponding self-map.

Summarizing, the ring of homotopy classes of natural transformations }
\[
\cancel{\Sigma^\infty B^{[2]} \to \Sigma^\infty B^{[2]}}
\] 
\st{isomorphic to the group ring $\Bbb Z[\Bbb Z_2]$.

With respect to this identification, an element $\theta \in \Bbb Z[\Bbb Z_2]$  corresponds to the homotopy class}
\begin{equation} \label{eqn:induced-theta}
\cancel{S[\Bbb Z_2]  \smsh_{\Bbb Z_2} B^{[2]} @> \hat{\theta} \smsh_{\Bbb Z_2} \text{id} >> S[\Bbb Z_2]  \smsh_{\Bbb Z_2} B^{[2]} \, .}
\end{equation}
\st{where $S[\Bbb Z_2]  \smsh_{\Bbb Z_2} B^{[2]} \simeq \Sigma^\infty B^{[2]}$. 

As a free  abelian group, $\Bbb Z[\Bbb Z_2]$ is generated by the elements
$1,\tau$, where $\tau$ is a generator of $\Bbb Z_2$.   
Let $\theta \in \Bbb Z[\Bbb Z_2]$ denote the element corresponding to $\lambda(2)$. Then}
\[
\cancel{\theta = a + b\tau}
\] 
\st{for suitable integers $a,b$.
Note that $\tau$ acts by switching the factors of $B^{[2]}$.

Let $W_{\theta}$ be the homotopy fiber of $\hat \theta$. Then $W_{\theta}$ is weakly contractible if and only if $\theta \in \Bbb Z[\Bbb Z_2]$ is a unit.

Consider the diagram}
\[
\cancel {\xymatrix{
B \ar[r] \ar@{=}[d] & \Omega \Sigma B \ar@{=}[d]\ar[r]^{\gamma} & Q(B^{[2]}) \ar[d]^{\theta} \\
B \ar[r]          & \Omega \Sigma B \ar[r]_{\lambda} & Q(B^{[2]})
} }
\]
\st{in which the rows are induces the metastable EHP sequences for $\gamma$ and $\lambda$ on homotopy groups.  The left right square of the diagram is homotopy commutative by definition of $\theta$.
Here we have used the same notation for the maps which induce the corresponding homotopy operations. The right vertical map is the map of infinite
loop spaces induced by  the map of spectra} \cancel{\eqref{eqn:induced-theta}.}

\st{The homotopy fiber of the right vertical map is the infinite loop space  $\Omega^\infty (W_{\theta}\smsh_{\Bbb Z_2} B^{[2]})$.
Hence the homogenous degree two functor $B \mapsto W_{\theta}\smsh_{\Bbb Z_2} B^{[2]}$ has vanishing homotopy groups in the metastable range:
If $B$ is $r$-connected, then $W_{\theta}\smsh_{\Bbb Z_2} B^{[2]}$ is $(3r+c)$-connected for some global constant $c$. 

If $B = S^{r+1}$
we have}
\[
\cancel{W_{\theta}\smsh_{\Bbb Z_2} (S^{r+1} \smsh S^{r+1}) = \Sigma^{r+1} W_{\theta} \smsh_{\Bbb Z_2} S^{(r+1)\alpha}}
\]
\st{implying that $W_{\theta} \smsh_{\Bbb Z_2} S^{(r+1)\alpha}$ is $2r+c$-connected (here $S^{(r+1)\alpha}$ is the one-point
compactification of $(r+1)$ copies of the sign representation). The group $\Bbb Z_2$ acts on $H_{r+1}(S^{(r+1)\alpha})$ by $(-1)^\epsilon$, where
$\epsilon = (-1)^{r+1}$. If we choose $r$ to be odd, then the first non-zero homology group of $\Sigma^{r+1} W_{\theta} \smsh_{\Bbb Z_2} S^{(r+1)\alpha}$
is in dimension $s+2r+3$, where $s$ is the connectivity of $W_\theta$. 
So we have $s+2r+3 \ge 3r+c$ or equivalently $s \ge r + c -3$. Since $r$ can vary,
we infer that  $s = \infty$. We infer that $W_\theta$ is weakly contractible.

Therefore $\theta\in \Bbb Z[\Bbb Z_2]$ is a unit. For  $\Bbb Z[\Bbb Z_2]$ the group units is isomorphic to Klein 4-group with elements}
\[
\cancel{\{\pm 1,\pm \tau\}\, .}
\]
\st{If $\theta \in \{-1,-\tau\}$, then the Cartan formula for $\lambda$ would be violated. To see this, let $1_B\: B\to B$ be the identity map. Then}
\[
\cancel{\lambda(1_B+1_B) = \lambda(1_B) + \Delta_B +\lambda(1_B) = \Delta_B}
\]
\st{where $\Delta_B\: B\to B\smsh B$ is the reduced diagonal. On the other hand}
\[
\cancel{\lambda(1_B+1_B)  = \theta\gamma(1_B + 1_B) = \theta \Delta_B}
\]
\st{So if $\theta = -1$,  we have $\Delta_B = \theta \Delta_B = -\Delta_B$. This identity is violated when $B = S^0$. Similarly,
if $\theta = -\tau$, we obtain $\Delta_B = \theta \Delta_B = -\tau\Delta_B = - \Delta_B$ and the identity is again violated
when $B = S^0$.}

\st{If $\theta = \tau$, then the Cartan formula for unstable maps $f,g \in [A,B]$ says}
\[
\cancel{\lambda (f+g) = f\cup g \, ,}
\]
\st{whereas if $\lambda = \tau \gamma$, one has}
\[
\cancel{\lambda(f+g) = \tau \gamma(f+g) = \tau (f\cup g) = g\cup f\, ,}
\]
\st{by the Cartan formula for $\gamma$.
Therefore $f\cup g = g\cup f$. 

However, the latter equation is violated in the following case: Let  $f,g\: S^1\times S^1\to S^1$ be the two projections.
Then  $f\cup g\: S^1\times S^1 \to S^1 \smsh S^1 = S^2$ has degree $+1$ and $g\cup f\: S^1\times S^1 \to S^1\smsh S^1 = S^2$ has
degree $-1$.

Therefore, $\theta\notin \{-1,\pm \tau\}$, so  $\theta = 1$. We conclude that  $\lambda = \gamma$.} \cancel{\qed }

{\color{blue} \begin{rem} Nick Kuhn has informed me
there are other papers in the literature that use
Goodwillie calculus ideas to give characterizations of James--Hopf
invariants. See e.g., \cite[App.~B]{Kuhn_diag}.

In Kuhn's scheme, a natural transformation  $\lambda\: J(B) \to Q(B \wedge B)$ is said to be a second James--Hopf invariant if the
 composition
\[
J_2(B) \hookrightarrow J(B) @> \lambda >> Q(B \smsh B)
\]
agrees with the composition
\[
J_2(B) \to B \smsh B @> \subset >> Q(B \smsh B).
\]
in the homotopy category of functors.

The ideas in the above proof show that $\lambda$ is completely characterized by this property.
\end{rem}}


\begin{thebibliography}{1}

\bibitem{Boardman-Steer}
J.~M. Boardman and B.~Steer, \emph{On {H}opf invariants}, Comment. Math. Helv.
  \textbf{42} (1967), 180--221.

\bibitem{Goodwillie_calc3}
Thomas~G. Goodwillie, \emph{Calculus. {III}. {T}aylor series}, Geom. Topol.
  \textbf{7} (2003), 645--711.

\bibitem{James_reduced_product}
I.~M. James, \emph{Reduced product spaces}, Ann. of Math. (2) \textbf{62}
  (1955), 170--197.

\bibitem{James-suspension-triad}
\bysame, \emph{On the suspension triad}, Ann. of Math. (2) \textbf{63} (1956),
  191--247.

\bibitem{klein2026stablehopf}
John~R. Klein, \emph{On the stable {H}opf invariant}, arXiv:2603.07854,
  submitted for publication.

\bibitem{Kuhn_diag}
Nicholas~J. Kuhn, \emph{Stable splittings and the diagonal}, Homotopy methods
  in algebraic topology ({B}oulder, {CO}, 1999), Contemp. Math., vol. 271,
  Amer. Math. Soc., Providence, RI, 2001, pp.~169--181.


\end{thebibliography}

\end{document}